\documentclass[12pt]{article}
\usepackage{amsmath, graphicx, amsfonts,amssymb, calrsfs}
\usepackage{color}

\def\bbR{\mathbb{R}}
\def\r{\rho}
\def\s{\sigma}

\def\dualorlicz{\widetilde{V}_{\phi}}
\def\vrad{\mathrm{vrad}}

\def\cK{\mathcal{K}}

\def\E{\mathcal{E}}

\def\o{\omega}
\def\ball{B^n_2}
\def\polar{K^\circ}
\def\l{\lambda}
\def\cS{\mathcal{S}}
\def\cL{\mathbf{L}}
\def\O{\Omega}
\def\bK{\mathbf{K}}
 
\def\DOrliczGmix{\widetilde{G}_{\vec{\phi}}^{orlicz}}
\def\DOrliczAmix{\widetilde{\O}_{\vec{\phi}}^{orlicz}}
\def\DOrliczGp{\widetilde{G}_{p}^{orlicz}}
\def\DOrliczAp{\widetilde{\O}_{p}^{orlicz}}
\def\DOrliczG{\widetilde{G}_{\phi}^{orlicz}}
\def\DOrliczA{\widetilde{\O}_{\phi}^{orlicz}}

\def\be{\begin{equation}}
\def\ee{\end{equation}}
\def\bea{\begin{eqnarray}}
\def\eea{\end{eqnarray}}
\def\bt{\begin{theorem}}
\def\et{\end{theorem}}
\def\bl{\begin{lemma}}
\def\el{\end{lemma}}
\def\br{\begin{remark}}
\def\er{\end{remark}}
\def\bc{\begin{corollary}}
\def\ec{\end{corollary}}
\def\bd{\begin{definition}}
\def\ed{\end{definition}}
\def\bp{\begin{proposition}}
\def\ep{\end{proposition}}
\newtheorem{theorem}{Theorem}[section]
\newtheorem{lemma}{Lemma}[section]
\newtheorem{remark}{Remark}[section]
\newtheorem{proposition}{Proposition}[section]
\newtheorem{corollary}{Corollary}[section]

\newtheorem{definition}{Definition}[section]
\begin{document}
\title{Dual Orlicz-Brunn-Minkowski theory: dual Orlicz $L_{\phi}$ affine and geominimal surface areas 
}

\author{Deping Ye }
\date{}
\maketitle
 
\begin{abstract} This paper aims to develop basic theory for the dual Orlicz $L_{\phi}$ affine and geominimal surface areas for star bodies, which belong to the recent dual Orlicz-Brunn-Minkowski theory for star bodies. Basic properties for these new affine invariants will be provided. Moreover,  related Orlicz affine isoperimetric inequality, cyclic inequality, Santal\'{o} style inequality and Alexander-Fenchel type inequality are established.  

 \end{abstract}\section{Introduction} \label{section introduction} 
The $L_p$ affine and geominimal surface areas are central in the $L_p$ Brunn-Minkowski theory for convex bodies (i.e., convex compact subsets of $\bbR^n$ with nonempty interiors). These affine invariants are very useful in applications, see \cite{Gr2, Jenkinson2012,   LudR, LR1, LSW, Paouris2010, SW5, Werner2012a, Werner2012b} among others.  Other major contributions, including the $L_p$ affine isoperimetric inequalities, can be found in, e.g.,  
  \cite{Lei1986, Lu1, MW1, MW2, Petty1974, Petty1985,  SW4, WY2008}.  Note that the $L_p$ affine and geominimal surface areas of $K$ for $p\geq 1$  in \cite{Lu1} were defined to  be (essentially)  the infimum of $V_p(K, L^\circ)$ with $L$ having the same volume as the unit Euclidean ball $\ball$ and with $L$ running over all star bodies and convex bodies respectively,  where $V_p(K, L^\circ)$ is the $p$-mixed volume of $K$ and the polar body of $L$. The author in  \cite{Ye2014a1} proved similar results for the $L_p$ affine surface area for $-n\neq p<1$, which motivate the definition of  the $L_p$ geominimal surface area for $-n\neq p<1$.  

There are dual concepts for the $L_p$ affine and geominimal surface areas, namely, the dual $L_p$ affine and geominimal surface areas for star bodies \cite{WH2008, WQ2011}, which belong to the dual ($L_p$) Brunn-Minkowski theory for star bodies developed by Lutwak \cite{Lut1975, Lut1988}. The dual ($L_p$) Brunn-Minkowski theory for star bodies received considerable attention, see  \cite{Bernig2014, GardnerDP2014, Gardner2007, GardnerV1998, GardnerV1999, GardnerV2000, Lutwak1990, Milman2006, Zhang1994} among others.  In particular, the dual ($L_p$) Brunn-Minkowski theory has been proved to be very powerful in solving many geometric problems, for instance, the Busemann-Petty problems (see e.g., \cite{Gardner1994, GardnerKoldobski1999, Lut1988, Zhang1999}). 
 
The Orlicz-Brunn-Minkowski theory for convex bodies, initiated from the work \cite{LYZ2010a, LYZ2010b} by Lutwak, Yang and Zhang, is the next generation of the $L_p$ Brunn-Minkowski theory for convex bodies. In view of the importance of the $L_p$ affine and geominimal surface areas in the $L_p$ Brunn-Minkowski theory, it is important to define Orlicz affine and geominimal surface areas.  Due to lack of homogeneity, extension of the $L_p$ affine and geominimal surface areas to their  Orlicz counterparts may not be unique. Here, we mention two major extensions in literature.  The first one is by Ludwig in \cite{Ludwig2009}, where  the general affine surface areas were proposed based on a beautiful integral expression of the $L_p$ affine surface areas. The second one is by the author in \cite{Ye2014}, where the Orlicz $L_{\phi}$ affine and geominimal surface areas were defined as the extreme values of $V_{\phi}(K, \vrad(L) L^\circ)$ with $L$ running over all star bodies and convex bodies  respectively. Readers are referred to \cite{Ludwig2009, Ye2013, Ye2014} for basic properties and inequalities regarding the Orlicz affine and geominimal surface areas.    
  
 This paper aims to develop the dual Orlicz $L_{\phi}$ affine and geominimal surface areas for star bodies, which belong to the recent dual Orlicz-Brunn-Minkowski theory for star bodies. Basic setting for  the dual Orlicz-Brunn-Minkowski theory has been developed in \cite{Ye2014a} (see the independent work \cite{ZZ} for special cases of some of the results in \cite{Ye2014a}), where the Orlicz radial addition was defined and the Orlicz $L_{\phi}$-dual mixed volume was proposed. Important inequalities in the classical Brunn-Minkowski theory, such as, Brunn-Minkowski inequality and Minkowski first inequality, have been extended to their dual Orlicz counterparts. In particular, the dual Orlicz-Minkowski inequality was proved and plays key roles in establishing Orlicz affine isoperimetric inequalities for the dual Orlicz $L_{\phi}$ affine and geominimal surface areas in this paper.

This paper is organized as follows.  Section \ref{Section: orlicz mixed volume} is dedicated to  the Orlicz $\phi$-mixed volume and its dual. In particular, the Orlicz isoperimetric inequality and the Orlicz-Urysohn inequality as well as their dual counterparts are proved.  In Section \ref{subsection dual affine}, the dual Orlicz $L_{\phi}$ affine and geominimal surface areas are proposed and their basic properties are proved. Related Orlicz affine isoperimetric inequality, Santal\'{o} style inequality and cyclic inequality are established. In Section \ref{section dual mixed},  the dual Orlicz mixed $L_{\phi}$ affine and geominimal surface areas for multiple star bodies are  briefly discussed and related Alexander-Fenchel type
inequality is given. Basic background and notation are provided in Section \ref{section background}, and more background can be found in  \cite{Gardner2006, Sch}.

% % % % % % % % % % % % % % % % % % % % % % % % % % % % % % % % % % % % % %
% ===================  Section 2 =========================================% 
% % % % % % % % % % % % % % % % % % % % % % % % % % % % % % % % % % % % % % 
 
\section{Background and Notation}\label{section background}
Denote by $\ball=\{x\in \bbR^n: \|x\|\leq 1\}$ the unit Euclidean ball in $\bbR^n$, where  $\|\cdot\|$ is the usual Euclidean metric induced by the inner product $\langle\cdot, \cdot\rangle$.  We use $\partial K$ to denote the boundary of $K$. In particular, $\partial \ball$ (usually denoted by  $S^{n-1}$) is the unit sphere in $\bbR^n$, and $S^{n-1}$ has the usual spherical measure $\s$.  In general, for a measurable set $K\subset \bbR^n$, $|K|$ denotes the Hausdorff content of the appropriate dimension of $K$.  For convenience, let $\o_n=|\ball|$.

The compact subset $K\subset \bbR^n$ is said to be star-shaped about the origin, if every closed line segment from the origin to any point $x\in K$ is contained in $K$.  Note that, if $K$ is star-shaped about the origin, then $K$ can be uniquely determined by its {\it radial function} $\r_K: S^{n-1}\rightarrow [0, \infty]$ defined by $ \r _K(u)=\max \{\l: \l u\in K\}$ for all $u\in S^{n-1}$. If $\r_K(u)$ is  continuous and positive on $S^{n-1}$, then  $K$ is said to be a star body (about the origin).  Let  $\cS_0$ denote the set of all star bodies (about the origin) in $\bbR^n$. Two star bodies $K, L\in \cS_0$ are dilates of each other if there is a constant $\lambda>0$ such that  $\r_L(u)=\lambda \r_K(u)$ for all $u\in S^{n-1}$;  equivalently,  $L=\lambda K=\{\lambda x,  x\in K\}.$ If $K\in \cS_0$ is convex, $K$ will be called a {\it convex body} (with the origin in its interior). The set of all convex bodies with the origin in their interiors is denoted by $\cK_0$ and clearly $\cK_0\subset \cS_0$.  Besides the radial function, a convex body $K\in \cK_0$ can be uniquely determined by its {\it support function}  $h_K(\cdot): S^{n-1}\rightarrow \bbR$ defined as $h_K(u)=\max _{x\in K} \langle x,u \rangle $ for all $ u\in S^{n-1}.$

 Define the {\it polar body}  $\polar$ of $K\in \cS_0$ by   $\polar=\{y\in \bbR^n: \langle x,y \rangle \leq 1, \forall x\in
K\}.$ It is easily checked that $\polar$ is always convex no matter whether $K\in \cS_0$ is convex or not. Note that $K\subset (K^\circ)^\circ$  for all $K\in \cS_0$. The bipolar theorem (see, e.g., \cite{Sch}) implies that, for $K\in \cS_0$, $(K^\circ)^\circ$ is equal to the convex hull of $K$ -- the smallest convex body contains $K$.  Moreover, if $K\in \cK_0$ is convex, $(K^\circ)^\circ =K$ and $\r_K(u){h_{K^\circ}(u)}=1$ holds for all $u\in S^{n-1}$.  

 Denote by $\cK_c$ and $\cK_s$ the sets of convex bodies with centroid and Santal\'{o} point at the origin, respectively.  Hereafter, $K\in \cK_0$ is said to have the {\it Santal\'{o} point} at the origin, if $\polar$ has the centroid at the origin, that is,  $K\in \cK_s \Leftrightarrow \polar \in \cK_c$.  For convenience, let $\widetilde{\cK}=\cK_c\cup \cK_s$ and  $\widetilde{\cS}=\{L\in \cS_0: \ L^\circ \in \widetilde{\cK} \}.$ Note that $K\in \widetilde{\cK}$ implies $K^\circ\in \widetilde{\cK}$.  Due to the bipolar theorem,  for $K\in \widetilde{\cS}$,  the convex hull of $K$ is a convex body in $\widetilde{\cK}$.  It is obvious that $\widetilde{\cK}\subset \widetilde{\cS}$. 
 
 For a linear transform $T: \bbR^n\rightarrow \bbR^n$,  $|det(T)|$, $T^{*}$ and $T^{-1}$ refer to the absolute value of the determinant, the transpose and the inverse of $T$ respectively. The set of all invertible linear transforms is denoted by $GL(n)$. 
We say $T\in SL(n)$ if $T\in GL(n)$ with $|det(T)|=1$. The set $T(K)$ with $K\in \cS_0$ will be written as $TK$ for simplicity.  An {\it origin-symmetric ellipsoid}  $\E\in \cK_0$ is the image of the Euclidean ball under some $T\in GL(n)$, that is, $\E=T\ball$ for some $T\in GL(n)$. Origin-symmetric ellipsoids serve as the maximizers/minimizers of many important affine isoperimetric inequalities, for example, the Blaschke-Santal\'{o} inequality:  for $K\in \widetilde{\cK}$,  $\vrad(K)\vrad(K^\circ)\leq 1$ with equality if and only if $K$ is an origin-symmetric ellipsoid. Hereafter, $\vrad(K)$ denotes the volume radius of $K$, i.e.,  \begin{equation*} \vrad(K)=\bigg(\frac{|K|}{|\ball|}\bigg)^{1/n}\Longleftrightarrow |K|^{1/n}=\o_n^{1/n}\vrad(K).\end{equation*} 
 Note that $\vrad(r\ball)=r$ for all $r>0$, and  for all $T\in SL(n)$,  \begin{equation}\label{affine:invariance:volume:radius} \vrad(TK)=\vrad(K).\end{equation}   It is easily checked that $|K|\leq |L|$ implies $\vrad(K)\leq \vrad(L)$. In particular,  $\vrad(K)\leq \vrad(L)$ if $K\subset L$. Due to $L\subset (L^\circ)^\circ$ for all $L\in \cS_0$, one gets: if $L\in \widetilde{\cS}$ (and hence $L^\circ\in \widetilde{\cK}$), then \begin{equation} \label{Mahler santalo-1}  \vrad(L)\vrad(L^\circ)\leq \vrad((L^\circ)^\circ)\vrad(L^\circ) \leq 1, \end{equation} with equality if and only if $L$ is an origin-symmetric ellipsoid. On the other hand, the following  inverse Santal\'{o} inequality \cite{BM} holds:   \begin{equation} \label{Mahler santalo-2}  \vrad(K)\vrad(K^\circ)\geq c, \ \ \ \ \forall K\in \widetilde{\cK}, \end{equation} where $c>0$  is a (universal) constant independent of $n$ and $K$  (see \cite{GK2, Nazarov2012} for estimates on $c$).

\section{Orlicz $\phi$-mixed volume and its dual} \label{Section: orlicz mixed volume}
       
       \subsection{Orlicz $\phi$-mixed volume}  Let $\phi: (0,\infty)\rightarrow (0, \infty)$ be a continuous function. 
  Define the Orlicz $\phi$-mixed volume  $V_{\phi}(K, Q)$ of convex bodies $K, Q\in \cK_0$ by \cite{Gardner2014, XJL, Ye2014}  $$V_{\phi}(K, Q)=\frac{1}{n}\int _{S^{n-1}}\phi\left(\frac{h_Q(u)}{h_K(u)}\right) h_K(u)\,dS(K, u),$$ where $S(K, \cdot)$ on $S^{n-1}$  is the surface area measure of $K$ (see \cite{Ale1937-1, Fenchel1938}), such that, for
any Borel subset $A$ of $S^{n-1}$, one has
$$S(K,A)=|\{x\in \partial K:\ \mbox{$\exists u\in A$, s.t., $H(x,u)$ is a support hyperplane of $\partial K$ at $x$}
\}|.$$    
 
  The following Orlicz-Minkowski inequality for $V_{\phi}(K, L)$ was established in \cite{Gardner2014} (where more general cases were also proved).  See \cite{XJL} for similar results. \bt \label{Minkowski inequality -1-1} {\bf (Orlicz-Minkowski inequality).} Let $K, L\in \cK_0$ and $\phi(t)$ be an  increasing convex function. One has,   $V_{\phi}(K, L)\geq |K|\cdot \phi\big({|L|^{1/n}}\cdot {|K|^{-1/n}}\big).$  If in addition $\phi(t)$ is strictly convex, equality holds if and only if $K$ and $L$ are dilates of each other.   \et 

 Define $ {S}_{\phi}(K)$, the  Orlicz $\phi$-surface area of $K$,  to be $nV_{\phi} (K, \ball)$. It is easily checked that for all $r>0$, ${S}_{\phi}(r \ball)= \phi(1/r)\cdot n|r\ball|.$ In particular,   \begin{equation*} {S}_{\phi}(B_K)= {S}_{\phi}\big(\vrad(K)\ball \big)= \phi\bigg(\frac{1}{\vrad(K)}\bigg)\cdot n|K|. \end{equation*}  The following result is an Orlicz isoperimetric inequality for ${S}_{\phi}(K)$, which follows from Theorem \ref{Minkowski inequality -1-1} by letting $L=\ball$.  The classical isoperimetric inequality is the special case with $\phi(t)=t$. 

\bt {\bf (Orlicz isoperimetric inequality).} Let $K\in \cK_0$ and $\phi(t)$ be an  increasing convex function. One has,  $S_{\phi}(K)\geq S_{\phi}(B_K).$  If $\phi(t)$ is strictly convex, equality holds if and only if $K$ is an origin-symmetric Euclidean ball.  \et

The following inequality is an Orlicz-Urysohn inequality for ${\o}_{\phi}(K)$, which follows from Theorem \ref{Minkowski inequality -1-1} by letting $K=\ball$ and $L=K$.  The classical Urysohn inequality is related to $\phi(t)=t$.  Here, ${\o}_{\phi}(K)$ is the  Orlicz ${\phi}$ mean width of $K\in \cK_0$  defined by  $${\o}_{\phi}(K)=\frac{1}{n\o_n} \int_{S^{n-1}}\phi\big({h_K(u)}\big)\,d\s(u)=\frac{V_{\phi}(\ball, K)}{\o_n}.$$  In particular, ${\o}_{\phi}(r\ball)=\phi(r)$ and hence ${\o}_{\phi}(B_K)=\phi(\vrad(K)).$  
\bt {\bf (Orlicz-Urysohn inequality).} Let $K\in \cK_0$ and $\phi(t)$ be an  increasing convex function. Then,   $\o_{\phi}(K)\geq \o_{\phi}(B_K).$  If $\phi(t)$ is strictly convex, equality holds if and only if $K$ is an origin-symmetric Euclidean ball.  \et

\subsection{Orlicz $L_\phi$-dual mixed volume} 
Let $\phi: (0,\infty)\rightarrow (0, \infty)$ be a continuous function.   
The Orlicz $L_{\phi}$-dual mixed volume for $K, L\in \cS_0$ was  defined in \cite{Ye2014a, ZZ}  by  $$\dualorlicz(K, L)=\frac{1}{n}\int _{S^{n-1}}\phi\left(\frac{\r_L(u)}{\r_K(u)}\right) [\r_K(u)]^n \,d\s(u).$$ Clearly, if $L=\lambda K$ for some $\lambda>0$, one gets \be \label{dual mixed L=aK} \dualorlicz (K, \lambda K)=\frac{1}{n}\int _{S^{n-1}}\phi\left( \lambda \right) [\r_K(u)]^n \,d\s(u)=\phi( \lambda) |K|.\ee  The following dual Orlicz-Minkowski inequality for $\dualorlicz(K, L)$  plays fundamental roles in this paper (more general results can be found in  \cite{Ye2014a}). See \cite{ZZ} for similar results.    
\bt \label{Minkowski inequality -1} {\bf (Dual Orlicz-Minkowski inequality).} Let $K, L\in \cS_0$, and  $F(t)=\phi(t^{1/n})$.  If $F(t)$ is concave, then $$\dualorlicz(K, L)\leq |K|\cdot \phi\bigg(\bigg(\frac{|L|} {|K|}\bigg)^{1/n}\bigg),$$ while if $F(t)$ is convex, the inequality is reversed. If  $F(t)$ is strictly concave (or convex, as appropriate), equality holds if and only if $K$ and $L$ are dilates of each other.   \et 
 
 Let  $\widetilde{S}_{\phi}(K)=n\dualorlicz (K, \ball)$ be the Orlicz $L_{\phi}$-dual surface area of $K$.  For $B_K=\vrad(K) \ball$ with $r=\vrad(K)$, one has, by formula (\ref{dual mixed L=aK}),  \begin{equation} \label{dual surface area of BK} \widetilde{S}_{\phi}(B_K)=\widetilde{S}_{\phi}(r \ball) = \phi(1/r)\cdot n|r\ball|=\phi\bigg(\frac{1}{\vrad(K)}\bigg)\cdot n|K|. \end{equation}  
 
 The following dual Orlicz isoperimetric inequality for $\widetilde{S}_{\phi}(K)$ follows immediately by letting $L=\ball$ in Theorem \ref{Minkowski inequality -1} and by formula (\ref{dual surface area of BK}).

\bt {\bf (Dual Orlicz isoperimetric inequality).} Let $K\in \cS_0$.   If  $F(t)=\phi(t^{1/n})$ is concave, then  $\widetilde{S}_{\phi}(K)\leq \widetilde{S}_{\phi}(B_K),$ while if $F(t)$ is convex, the inequality is reversed. If  $F(t)$ is strictly concave (or convex, as appropriate),  equality holds if and only if $K$ is an origin-symmetric Euclidean ball. \et
 
    \noindent {\bf Remark.} Write $\widetilde{S}_p(K)$ for the case  $\phi(t)=t^p$, and $\widetilde{S}_p(\lambda K)=\lambda^{n-p}\widetilde{S}_p(K),$ for all $\lambda>0$.  When $\phi(t)=t^p$ with $p\in [n, \infty)\cup (-\infty, 0)$, one can even have $$\frac{\widetilde{S}_p(K)}{\widetilde{S}_p(\ball)}\geq \bigg(\frac{|K|}{|\ball|}\bigg)^{\frac{n-p}{n}}, $$ while if $p\in (0, n)$, the inequality is reversed.  
 
For $K\in \cS_0$, define the Orlicz $L_{\phi}$ mean radius of $K$ (denoted by $\widetilde{\o}_{\phi}(K)$)  as $$\widetilde{\o}_{\phi}(K)=\frac{1}{n\o_n} \int_{S^{n-1}}\phi\big ( {\r_K(u)}\big)\,d\s(u)=\frac{\dualorlicz(\ball, K)}{\o_n}.$$  The following dual Orlicz-Urysohn inequality for $\widetilde{\o}_{\phi}(K)$ follows immediately by letting $K=\ball$ and $L=K$ in Theorem \ref{Minkowski inequality -1}.  
\bt {\bf (Dual Orlicz-Urysohn inequality).} Let $K\in \cS_0$.  If  $F(t)=\phi(t^{1/n})$ is concave, then $\widetilde{\o}_{\phi}(K)\leq \widetilde{\o}_{\phi}(B_K),$  while if $F(t)$ is convex, the inequality is reversed. If  $F(t)$ is strictly concave (or convex, as appropriate),  equality holds if and only if $K$ is an origin-symmetric Euclidean ball. \et

% % % % % % % % % % % % % % % % % % % % % % % % % % % % % % % % % % % % % %
% ===================  Section 3 =========================================% 
% % % % % % % % % % % % % % % % % % % % % % % % % % % % % % % % % % % % % % 

\section{Dual Orlicz $L_{\phi}$ affine and geominimal surface areas}\label{subsection dual affine}
\subsection{Definitions and basic properties} 
Let $\phi: (0, \infty)\rightarrow (0, \infty)$ be a continuous function. Consider the sets of functions $\widetilde{\Phi}$ and $\widetilde{\Psi}$ \begin{eqnarray*} \widetilde{\Phi}&=&\{\phi :  \mbox{$F(t)=\phi(t^{1/n})$ is either a constant or a strictly convex function}\}, \\  \widetilde{\Psi}&=&\{\phi:   \mbox{$F(t)=\phi(t^{1/n})$ is either a constant or an increasing strictly concave function}\}. \end{eqnarray*}   Note that $t^p$ with $p\in (-\infty, 0)\cup (n, \infty)$ and  all decreasing  strictly convex functions are in $\widetilde{\Phi}$; while  $t^p$ with $p\in (0, n)$ and all increasing strictly concave functions are in $\widetilde{\Psi}$.

For $K\in \cS_0$, denote by $\DOrliczA(K)$  the dual Orlicz $L_{\phi}$ affine surface area of $K$ and  by $\DOrliczG(K)$ the dual Orlicz $L_{\phi}$ geominimal surface area of $K$.  

  \bd \label{Orlicz affine surface}  Let $K\in \cS_0$ be a star  body about the origin. \vskip 2mm \noindent
 (i)  For $\phi\in \widetilde{\Phi}$, define $\DOrliczA(K)$ and $\DOrliczG(K)$ as follows: 
\begin{eqnarray}
  \DOrliczA(K)=\inf_{L\in \widetilde{\cS}} \left\{n\dualorlicz(K, \vrad(L^\circ)L ) \right\}, \ \ \& \ \ 
    \DOrliczG(K)=\inf_{L\in \widetilde{\cK} } \left\{n\dualorlicz(K, \vrad(L^\circ)L ) \right\}. \label{Orlicz-geominimal-surface-area-convex}
 \end{eqnarray}  
  (ii)  For $\phi\in \widetilde{\Psi}$, define $\DOrliczA(K)$ and $\DOrliczG(K)$ as follows:  
  \begin{eqnarray}
  \DOrliczA(K)=\sup_{L\in \widetilde{\cS}} \left\{n\dualorlicz(K, \vrad(L^\circ)L ) \right\}, \ \ \& \ \   \DOrliczG(K)=\sup_{L\in \widetilde{\cK}} \left\{n\dualorlicz(K, \vrad(L^\circ)L ) \right\}. \label{Orlicz-geominimal-surface-area:concave} 
  \end{eqnarray}  \ed

\noindent {\bf Remark.}  Note that if $\phi(t)=\alpha$ is a constant function, then $\DOrliczA(K)=\DOrliczG(K)=\alpha\cdot n|K|$ for all $K\in \cS_0$.  It is often more convenient to take the infimum/supremum over $L$ with $|L^\circ|=\o_n$. In fact, for all $L\in \cS_0$, one checks that  $|\big(\vrad(L^\circ)L\big)^\circ|=\big|\frac{1}{\vrad(L^\circ)}L^\circ\big|=\o_n$. Hence, for $\phi\in \widetilde{\Phi}$,  \begin{eqnarray*} 
    \DOrliczG(K)=\inf_{L\in \widetilde{\cK} } \left\{n\dualorlicz(K, \vrad(L^\circ)L ) \right\}=\inf \left\{n\dualorlicz(K, L):  L\in \widetilde{\cK} \ with\ |L^\circ|=\o_n \right\}.  
 \end{eqnarray*}  Similar formulas for other cases can be obtained along the same line.  It is easy to prove that, for all $K\in \cS_0$ and  $\phi\leq \psi$ with either $\phi, \psi\in \widetilde{\Phi}$   or  $\phi, \psi\in \widetilde{\Psi}$, then $$
\DOrliczA(K)\leq  \widetilde{\O}_{\psi}^{orlicz} (K), \ \ \& \ \  
\DOrliczG(K)\leq  \widetilde{G}_{\psi}^{orlicz} (K).$$ Moreover, by $\widetilde{\cK}\subset \widetilde{\cS}$ and by taking $L=\ball$ in Definition \ref{Orlicz affine surface}, one has, for all $K\in \cS_0$,  
\begin{eqnarray*}  \DOrliczA(K) &\leq& \DOrliczG(K)\leq \widetilde{S}_{\phi}(K), \ \ \ \forall \phi\in \widetilde{\Phi};\\    \DOrliczA(K)&\geq& \DOrliczG(K)\geq \widetilde{S}_{\phi}(K),  \ \ \ \forall \phi\in \widetilde{\Psi}.\end{eqnarray*}

 \vskip 2mm 
We now prove that the dual Orlicz $L_{\phi}$ affine and geominimal surface areas are affine invariant. 
 \bp\label{homogeneous:degree}  Let $K\in \cS_0$. For all $\phi\in \widetilde{\Phi}\cup \widetilde{\Psi}$, one has $$ \DOrliczA(TK)=\DOrliczA(K); \ \ \ \DOrliczG(TK)=\DOrliczG(K), \ \  \forall T\in SL(n).$$ \ep 
\noindent {\bf Proof.}    Let $T\in SL(n)$ and let $u=\frac{T^{-1} v}{\|T^{-1} v\|}\in S^{n-1}$ for $v\in S^{n-1}$.   It can be checked that $\r_{TK}(v)\|T^{-1}v\|=\r_K(u)$ (see \cite{Gardner2006}), and $ [\r_{TK}(v)]^n\,d\s(v)=[\r_K(u)]^n\,d\s(u)$. Hence   \begin{eqnarray*}
n\dualorlicz(TK, TL)&=& \int _{S^{n-1}}  \phi \bigg( \frac{\r_{TL}(v)}{\r_{TK}(v)}\bigg) [\r_{TK}(v)]^n\,d\s(v)\\ &=&   \int _{S^{n-1}}  \phi\bigg( \frac{\r_{L}(u)}{\r_{K}(u)}\bigg) [\r_{K}(u)]^n\,d\s(u) =n  \dualorlicz(K, L). \end{eqnarray*} 
 Together with formula (\ref{affine:invariance:volume:radius}), one gets, $\forall \phi\in \widetilde{\Phi},$ \begin{eqnarray*}
\DOrliczG(TK) &=& \inf_{L\in \widetilde{\cK}} \left\{ n\dualorlicz\big(TK, \vrad((TL)^\circ)(TL)\big)  \right\}\\ &=&\inf_{L\in \widetilde{\cK}} \left\{ n\dualorlicz(K, \vrad(L^\circ) L)  \right\}=\DOrliczG (K).
\end{eqnarray*} The remainder of the theorem follows along the same line.

\vskip 2mm \noindent {\bf Remark.} Let $\DOrliczAp(K)=\DOrliczA(K)$ and $\DOrliczGp(K)=\DOrliczG(K)$ for $\phi(t)=t^p$ with $n\neq p\in \bbR$. A more careful calculation shows that  $$ \DOrliczGp (TK)=|det(T)|^{\frac{n-p}{n}} \DOrliczGp(K); \ \ \DOrliczAp (TK)=|det(T)|^{\frac{n-p}{n}} \DOrliczAp(K).$$

 Denote by $\widetilde{\Phi}_1$ the set of functions $\phi(t)\in \widetilde{\Phi}$ with $F(t)=\phi(t^{1/n})$ being either a constant function or a decreasing strictly convex function. Clearly, $\phi\in \widetilde{\Phi}_1$ is decreasing. 
 \bc\label{Orlicz:ellpsoids} Let $\E$ be an origin-symmetric ellipsoid. For $\phi\in\widetilde{\Phi}_1\cup \widetilde{\Psi}$, one has   $$\DOrliczA(\E)=\DOrliczG(\E)=\phi\left( \frac{1} {\vrad(\E )} \right) \cdot n |\E|. $$ In particular,  $ \DOrliczG( \E)= \DOrliczA( \E)=  n| \E| \cdot  \phi\left( 1 \right),$ if $|\E|=|\ball|$. 
\ec   \noindent {\bf Proof.} For $\E=T(r\ball)$ with $T\in SL(n)$ and $r>0$,  Proposition \ref{homogeneous:degree} implies that, for $\phi\in \widetilde{\Phi}_1\cup \widetilde{\Psi}$,  
  $$\DOrliczA(\E) = \DOrliczA(r\ball), \ \ \&\ \ \DOrliczG(\E)=\DOrliczG(r\ball).$$   
   It is enough to consider $\E=r\ball$ for some $r>0$. By formulas (\ref{Orlicz-geominimal-surface-area-convex}) and (\ref{dual surface area of BK}), one has, for $\phi\in \widetilde{\Phi}_1$, \begin{eqnarray*}
  \DOrliczA(r\ball) \leq  \DOrliczG(r\ball) \leq  n  \dualorlicz(r\ball ,  \ball)  =n \phi (1/r) |r\ball|.
  \end{eqnarray*} On the other hand,  dual Orlicz-Minkowski inequality and inequality (\ref{Mahler santalo-1}) together with the decreasing property of $\phi(t)$ (as $F(t)$ is decreasing) imply that   \begin{eqnarray*}
     \DOrliczA(r\ball) &=& n \inf_{L\in \widetilde{\cS}}\dualorlicz(r\ball , \vrad(L^\circ) L) \\ & \geq&    \inf_{L\in \widetilde{\cS}} \phi\left(\! \frac{\vrad(L)\vrad(L^\circ)} {r}\!\right)\cdot n |r\ball|  \geq   n \phi\left( 1/ r\right) |r\ball|.\end{eqnarray*}  Hence,  $
     \DOrliczG(r\ball)= \DOrliczA(r\ball)=  n  \phi\left(1/r\right) |r\ball|$  for $\phi\in\widetilde{\Phi}_1$ as desired.  
     
      Similarly,  dual Orlicz-Minkowski inequality and inequality (\ref{Mahler santalo-1})   imply that for $\phi\in\widetilde{\Psi}$,  \begin{eqnarray*}
n \phi (1/r) |r\ball| & \leq &   \DOrliczG(r\ball)\leq \DOrliczA(r\ball)  \\ &\leq &     \sup_{L\in \widetilde{\cS}} \phi\left( \frac{\vrad(L)\vrad(L^\circ)}{r} \right)\cdot n |r\ball|  \leq    n \phi\left(1/ r\right) |r\ball|.\end{eqnarray*}      

\subsection{Orlicz affine isoperimetric inequalities} \label{subsection inequality} 
 
\bp\label{bounded by volume}  Let $K\in \cS_0$ be a star body about the origin and let $\phi\in \widetilde{\Phi}$. Then \begin{eqnarray*}   \DOrliczA(K) &\leq&    \phi\big({vrad(K^\circ)}\big) \cdot n|K|, \ \  for \ \ K\in \widetilde{\cS};\\  
 \DOrliczA(K) &\leq&   \DOrliczG(K) \leq   \phi\big({vrad(K^\circ)}\big) \cdot n|K|, \ \ for \ \ K\in \widetilde{\cK}. \end{eqnarray*} 
 Moreover,  if  in addition $\phi\in  \widetilde{\Phi}_1$, one has, for $K\in \cS_0$,  $$\DOrliczG(K)\geq \DOrliczA(K) \geq    \phi\bigg(\frac{1}{\vrad(K)}\bigg) \cdot n|K|.$$ If $\phi\in  \widetilde{\Psi}$, the above inequalities are reversed.    
\ep

\noindent {\bf Proof.}  We only prove the case $\phi\in \widetilde{\Phi}$, and the case $\phi\in \widetilde{\Psi}$ follows along the same fashion.  Formula (\ref{dual mixed L=aK}) implies that,  \begin{eqnarray*} \DOrliczA(K) \!\!\!&=&\!\!\!\inf_{L\in \widetilde{\cS} } \left\{n\dualorlicz(K, \vrad(L^\circ)L ) \right\}\leq   n\dualorlicz(K, \vrad(K^\circ)K )   =\phi\left( {\vrad(K^\circ)}\right)\cdot n |K|,  \  \ K\in \widetilde{\cS}, \\  \DOrliczA(K)\! \!\!&\leq&\!\!\! \DOrliczG(K) \leq   n\dualorlicz(K, \vrad(K^\circ)K ) =\phi\left( {\vrad(K^\circ)}\right)\cdot n |K|,  \ \ \ K\in \widetilde{\cK}. \end{eqnarray*}   
 Now assume $\phi\in \widetilde{\Phi}_1$, and hence $\phi(t)$ is decreasing.  Together with dual Orlicz-Minkowski inequality and inequality (\ref{Mahler santalo-1}), one has, for all $L\in \widetilde{\cS}$,  \begin{eqnarray*}  \DOrliczG(K) &\geq& \DOrliczA(K) =\inf_{L\in \widetilde{\cS}} \{n\dualorlicz(K, \vrad(L^\circ) L)\} \\ &\geq&   n|K|\cdot  \inf_{L\in \widetilde{\cS}}  \phi\left( \frac{\vrad(L)\vrad(L^\circ)} { \vrad(K)}   \right) \geq   \phi\bigg(\frac{1}{\vrad(K)}\bigg) \cdot n|K|.   \end{eqnarray*}

We prove the following Orlicz affine isoperimetric inequalities for the dual Orlicz $L_{\phi}$ affine and geominimal surface areas. Let $B_K=\vrad(K)\ball$.  Corollary \ref{Orlicz:ellpsoids} implies that, for $\phi\in  \widetilde{\Phi}_1\cup  \widetilde{\Psi}$,
\begin{equation} \label{equal to the ball}
 \phi\bigg(\frac{1}{\vrad(K)}\bigg) \cdot n|K|=\phi\bigg(\frac{1}{\vrad(B_K)}\bigg) \cdot n|B_K| =  \DOrliczA(B_K)=\DOrliczG(B_K).\end{equation}

 \bt \label{isoperimetric:geominimal} Let $K\in \cS_0$ be a star body about the origin. 
\vskip 2mm \noindent (i)  For $\phi\in \widetilde{\Phi}_1$ and $K\in \cS_0$, one has $$\DOrliczG(K)\geq \DOrliczA(K) \geq    \DOrliczA(B_K)=\DOrliczG(B_K), $$ while if  $\phi\in \widetilde{\Psi}$, the inequalities are reversed. 

 \noindent 
(ii)  For $\phi\in \widetilde{\Phi}_1$, one has \begin{eqnarray*} \DOrliczA(K)&\leq& \DOrliczA((B_{K^\circ})^\circ),\ \ \ for \ K\in \widetilde{\cS}; \\ 
 \DOrliczA(K)&\leq&  \DOrliczG(K) \leq    \DOrliczG((B_{K^\circ})^\circ),  \ \ for \ K\in \widetilde{\cK}.\end{eqnarray*}
Equality holds if and only if $K$ is an origin-symmetric ellipsoid.\vskip 2mm \noindent 
 \et

\noindent {\bf Proof.} 
(i). Let $\phi\in \widetilde{\Phi}_1$.  For all $K\in \cS_0$, Proposition \ref{bounded by volume} and equality (\ref{equal to the ball}) imply that 
 $$\DOrliczG(K)\geq \DOrliczA(K) \geq    \phi\bigg(\frac{1}{\vrad(K)}\bigg) \cdot n|K| =\DOrliczA(B_K)=\DOrliczG(B_K).$$  The case  $\phi\in \widetilde{\Psi}$ follows along the same fashion.

\vskip 2mm \noindent (ii). Let $\phi\in \widetilde{\Phi}_1$.  Proposition \ref{bounded by volume}, Corollary \ref{Orlicz:ellpsoids} and inequality (\ref{Mahler santalo-1}) imply that for all $K\in \widetilde{\cS}$,   \begin{eqnarray} \DOrliczA(K)&\leq& \phi\left( {\vrad(K^\circ)}\right)\cdot n |K|=\phi\left({\vrad(B_{K^\circ})}\right)\cdot n |(B_{K^\circ})^\circ|\cdot \frac{|K|}{ |(B_{K^\circ})^\circ|}\nonumber\\ &=&\phi\bigg(\frac{1}{\vrad\big((B_{K^\circ})^\circ\big)}\bigg)\cdot n |(B_{K^\circ})^\circ|\cdot \frac{|K||K^\circ|}{ |(B_{K^\circ})^\circ| |B_{K^\circ}|}\leq \DOrliczA((B_{K^\circ})^\circ).\label{Phi-1 isoperimetric} \end{eqnarray}   Clearly, equality holds if $K$ is an origin-symmetric ellipsoid. On the other hand, to have equality in the above inequalities, one needs to have equality in the second inequality. That is, equality holds in inequality (\ref{Mahler santalo-1}), and hence $K$ has to be an origin-symmetric ellipsoid.  Similarly, for $K\in \widetilde{\cK}$, one has,  $ \DOrliczG(K)  \leq  \phi\left({\vrad(K^\circ)}\right)\cdot n |K|  \leq  \DOrliczG((B_{K^\circ})^\circ),$  with equality if and only if $K$ is an origin-symmetric ellipsoid.

  \vskip 2mm \noindent {\bf Remark.}  Part (i) of Theorem \ref{isoperimetric:geominimal} asserts that {\em among all $K\in \cS_0$ with volume fixed, the dual Orlicz $L_{\phi}$ affine and geominimal surface areas attain the minimum  for $\phi\in \widetilde{\Phi}_1$ and the maximum for  $\phi\in \widetilde{\Psi}$ at origin-symmetric ellipsoids.}     For $\phi\in \widetilde{\Psi}$ and for $K\in \widetilde{\cK}$, one can prove that, \begin{eqnarray*}\DOrliczA(K) \geq  \DOrliczG(K) \geq n |(B_{K^\circ})^\circ|\cdot \phi\bigg(\frac{1}{\vrad\big((B_{K^\circ})^\circ\big)}\bigg)  \cdot  \frac{|K||K^\circ|}{ |(B_{K^\circ})^\circ| |B_{K^\circ}|} \geq   c^n \cdot \DOrliczG\big((B_{K^\circ})^\circ\big), \end{eqnarray*} where $c$ is the constant in inequality (\ref{Mahler santalo-2}).  Similar to inequality (\ref{Phi-1 isoperimetric}), for all $\phi\in \widetilde{\Phi}\setminus \widetilde{\Phi}_1$, \begin{eqnarray*} \DOrliczA(K) &\leq&   \phi\bigg(\frac{1}{\vrad\big((B_{K^\circ})^\circ\big)}\bigg)\cdot n |(B_{K^\circ})^\circ|, \ \ \ K\in \widetilde{\cS}, \\ 
  \DOrliczG(K) &\leq&   \phi\bigg(\frac{1}{\vrad\big((B_{K^\circ})^\circ\big)}\bigg)\cdot n |(B_{K^\circ})^\circ|, \ \ \ K\in \widetilde{\cK}. \end{eqnarray*}

   \subsection{Santal\'{o} style inequalities} 
 
 The following proposition gives Santal\'{o} style inequalities for   $\DOrliczAp(K)$ and  $\DOrliczGp(K)$.  
  \bp 
Let $n\neq p\in \bbR$. 
\vskip 2mm \noindent 
  (i) Let $-n\leq p\leq 0$. For $K\in \widetilde{\cK}$, one has  \begin{eqnarray*} &&c^{n-p} \big[ \DOrliczAp(\ball) \big]^2  \leq   \DOrliczAp(K)  \DOrliczAp(K^\circ)\leq   \DOrliczGp(K)  \DOrliczGp(K^\circ)\leq \big[\DOrliczGp(\ball)\big]^2, \end{eqnarray*}  with equality in $\DOrliczGp(K)  \DOrliczGp(K^\circ)\leq \big[\DOrliczGp(\ball)\big]^2$ if and only if $K$ is an origin-symmetric ellipsoid.   Moreover, $\DOrliczAp(K)  \DOrliczAp(K^\circ)  \leq \big[\DOrliczAp(\ball)\big]^2$ for $K\in \widetilde{\cS}$  with equality if and only if $K$ is an origin-symmetric ellipsoid. 
  
\vskip 2mm \noindent  (ii)  Let $p<- n$. For $K\in \widetilde{\cK}$, one has  \begin{eqnarray*} &&c^{n-p}  \big[\DOrliczAp(\ball)\big]^2  \leq   \DOrliczAp(K)  \DOrliczAp(K^\circ)\leq   \DOrliczGp(K)  \DOrliczGp(K^\circ)\leq  c^{n+p} \big[\DOrliczGp(\ball)\big]^2. \end{eqnarray*}      
\noindent (iii) Let $0< p<n$. For $K\in \widetilde{\cK}$, one has  \begin{eqnarray*} &&c^{n+p} \big[ \DOrliczGp(\ball) \big]^2  \leq   \DOrliczGp(K)  \DOrliczGp(K^\circ)\leq   \DOrliczAp(K)  \DOrliczAp(K^\circ)\leq  \big[\DOrliczAp(\ball)\big]^2. \end{eqnarray*}    Furthermore, $ \DOrliczAp(K)  \DOrliczAp(K^\circ)  \leq \big[\DOrliczAp(\ball)\big]^2 $  for $K\in \widetilde{\cS}$   with equality if and only if $K$ is an origin-symmetric ellipsoid.

\vskip 2mm \noindent (iv) Let $p>n$. The following inequalities hold \begin{eqnarray*}   \DOrliczGp(K)  \DOrliczGp(K^\circ) &\leq&  (n\o_n)^2, \ \ \ K\in \widetilde{\cK}, \\     \DOrliczAp(K)  \DOrliczAp(K^\circ)  &\leq& (n\o_n)^2, \ \ \  K\in \widetilde{\cS}. \end{eqnarray*}     \ep

   \noindent {\bf Proof.}  Replacing $K\in \widetilde{\cK}$ by its polar body $K^\circ\in \widetilde{\cK}$ in  Proposition \ref{bounded by volume}, one has, for $\phi\in \widetilde{\Phi}$,
$$  \DOrliczA(K^\circ) \leq  \DOrliczG(K^\circ) \leq \phi\!\left( {\vrad(K)}\!\right) \cdot n|K^\circ|.$$ Hence, for $\phi\in \widetilde{\Phi}$ and $K\in \widetilde{\cK}$, \begin{eqnarray}   
  \DOrliczA(K)\DOrliczA(K^\circ)  \leq  \DOrliczG(K)\DOrliczG(K^\circ)\leq  \phi\!\left( {\vrad(K)}\!\right)\phi\!\left( {\vrad(K^\circ)}\!\right) \cdot n^2|K|\cdot |K^\circ|.\ \ \ \  \label{Mahler:product--1}
\end{eqnarray}   
Moreover, for $\phi\in \widetilde{\Phi}_1$ and $K\in \cS_0$, 
 \begin{eqnarray}
  \phi\!\left(\!\frac{1} {\vrad(K)}\!\right)\phi\!\left(\! \frac{1}{\vrad(K^\circ)}\!\right) \cdot n^2|K|\cdot |K^\circ|  \leq  \DOrliczA(K)\DOrliczA(K^\circ)  \leq  \DOrliczG(K)\DOrliczG(K^\circ). \ \ \ \  \label{Mahler:product--3}
\end{eqnarray}

\noindent (i).  For $-n\leq p\leq 0$, one gets $n-p\geq n+p\geq 0$ and $\phi(t)=t^p\in \widetilde{\Phi}_1$.   Inequalities (\ref{Mahler:product--1}) and (\ref{Mahler:product--3})  together with Corollary \ref{Orlicz:ellpsoids}, the Blaschke-Santal\'{o} inequality and the inverse Santal\'{o} inequality imply that for $K\in \widetilde{\cK}$,  \begin{eqnarray*} c^{n-p}\cdot \big[\DOrliczAp(\ball)\big]^2 \!\!\!&=&\!\!\!   c^{n-p}\cdot n^2 |\ball|^2 \leq  \frac{n^2(|K|\cdot |K^\circ|)^{\frac{n-p}{n}}}{|\ball|^{\frac{-2p}{n}}}  \leq \DOrliczAp(K)\DOrliczAp(K^\circ)  
    \\ \!\!\!&\leq&\!\!\!  \DOrliczGp(K)\DOrliczGp(K^\circ) \leq \frac{n^2(|K|\cdot |K^\circ|)^{\frac{n+p}{n}}}{|\ball|^{\frac{2p}{n}}}  \leq n^2 |\ball|^2=\big[\DOrliczGp(\ball)\big]^2.
\end{eqnarray*}  The equality clearly holds in $ \DOrliczGp(K)  \DOrliczGp(K^\circ)\leq  \big[\DOrliczGp(\ball)\big]^2$ if $K$ is an origin-symmetric ellipsoid.   On the other hand, equality holds  only if equality holds in the Blaschke-Santal\'{o} inequality, that is,   $K$ has to be an origin-symmetric ellipsoid.
 
The proof of $\DOrliczAp(K)  \DOrliczAp(K^\circ)\leq  \big[\DOrliczAp(\ball)\big]^2$ for $K\in \widetilde{\cS}$ with characterization for equality follows along the same line and hence is omitted. 

\vskip 2mm \noindent (ii). For $p<-n$, one gets $n+p<0<n-p$ and $\phi(t)\in \widetilde{\Phi}_1$.  Similar to part (i),   Corollary \ref{Orlicz:ellpsoids} and the inverse Santal\'{o} inequality imply that for all $K\in \widetilde{\cK}$,  \begin{eqnarray}
    c^{n-p}\cdot \big[\DOrliczAp(\ball)\big]^2\!\!\!&\leq &\!\!\!   \frac{n^2(|K|\cdot |K^\circ|)^{\frac{n-p}{n}}}{|\ball|^{-\frac{2p}{n}}} \leq \DOrliczAp(K)\DOrliczAp(K^\circ) \nonumber \\  \!\!\!&\leq&\!\!\!    \DOrliczGp(K)\DOrliczGp(K^\circ)  \leq   \frac{n^2(|K|\cdot |K^\circ|)^{\frac{n+p}{n}}}{|\ball|^{\frac{2p}{n}}} \leq  c^{n+p}\cdot \big[\DOrliczGp(\ball)\big]^2. \nonumber \end{eqnarray}  
  \noindent  (iii). Let $0< p< n$ which implies $n+p>n-p>0$ and  $\phi(t)=t^p\in\widetilde{\Psi}$.    By Corollary \ref{Orlicz:ellpsoids}, the Blaschke-Santal\'{o} and the inverse Santal\'{o} inequalities, one gets, for $K\in \widetilde{\cK}$,  \begin{eqnarray*} \big[\DOrliczAp(\ball)\big]^2 &\geq& \frac{n^2(|K|\cdot |K^\circ|)^{\frac{n-p}{n}}}{|\ball|^{-\frac{2p}{n}}}  \geq \DOrliczAp(K)\DOrliczAp(K^\circ)  
    \\ &\geq&  \DOrliczGp(K)\DOrliczGp(K^\circ) \geq \frac{n^2(|K|\cdot |K^\circ|)^{\frac{n+p}{n}}}{|\ball|^{\frac{2p}{n}}}  \geq     c^{n+p} \big[\DOrliczGp(\ball)\big]^2.
\end{eqnarray*}  
Similarly,  $ \DOrliczAp(K)  \DOrliczAp(K^\circ)\leq  \big[\DOrliczAp(\ball)\big]^2$ for $K\in \widetilde{\cS}$. Moreover, equality holds only if equality holds in inequality (\ref{Mahler santalo-1}) and hence $K$ has to be an origin-symmetric ellipsoid. On the other hand, equality clearly holds for all origin-symmetric ellipsoids.   
  
 \vskip 2mm \noindent (iv).  Let $p>n$. Similar to part (i), one has,  by inequality (\ref{Mahler santalo-1}), 
 \begin{eqnarray*}
    \DOrliczAp(K)\DOrliczAp(K^\circ) & \leq & \frac{n^2(|K|\cdot |K^\circ|)^{\frac{n+p}{n}}}{|\ball|^{\frac{2p}{n}}}\leq  (n\o_n)^2, \ \ \ \ K\in \widetilde{\cS}.  \end{eqnarray*}      Similarly, $  \DOrliczGp(K)\DOrliczGp(K^\circ)   \leq  (n\o_n)^2$ for $K\in \widetilde{\cK}$.

    \subsection{Cyclic inequalities and a monotonicity property} \label{subsection cyclic inequalities}  
     Let $H(t)= (\phi\circ \psi^{-1})(t)$ be the composition of $\phi(t)$ and $\psi^{-1}(t)$, where $\psi^{-1}(t)$, the inverse function of $\psi(t)$, is always assumed to exist. Let $H(0)=\lim_{t\rightarrow 0} H(t)$ if the limit exists and is finite; while let $H(0)=\infty$ if $\lim_{t\rightarrow 0} H(t)=\infty$. Similarly, let $H(\infty)=\lim_{t\rightarrow \infty} H(t)$ if the limit exists and is finite; or simply $H(\infty)=\infty$ if $\lim_{t\rightarrow \infty} H(t)=\infty$.  As  in \cite{Ye2014}, we are not interested in the following cases: $H(t)$ being decreasing with $\phi(t), \psi(t)\in \widetilde{\Psi}$ (as all functions  $\phi(t)\in\widetilde{\Psi}$ are increasing and hence $H(t)$ is always increasing), and  $H(t)$ being concave decreasing (as otherwise $\phi$ is eventually a constant function). 
  
\bt\label{cyclic}
Let $K\in \cS_0$ and $H(t)$ be as above. 
\vskip 2mm \noindent
(i) Assume that $\phi$ and $ \psi$ satisfy one of the following conditions: (a)  $\phi\in \widetilde{\Phi}$ and $\psi\in \widetilde{\Psi}$ with $H(t)$ increasing; (b) $\phi, \psi\in \widetilde{\Phi}$ with $H(t)$ decreasing. Then, \begin{eqnarray*}
\frac{\DOrliczA(K)}{n|K|}  \leq   H\bigg(\frac{\widetilde{\O}_{\psi}^{orlicz}(K)}{n|K|} \bigg)   \  if \ K\in \widetilde{\cS}, \ \ and\   \  \frac{\DOrliczG(K)}{n|K|}   \leq   H\bigg(\frac{\widetilde{G}_{\psi}^{orlicz}(K)}{n|K|} \bigg)  \  if \ K\in \widetilde{\cK}. \end{eqnarray*} While if $\phi$ and $ \psi$ satisfy condition (c) $\phi\in \widetilde{\Psi}$ and $\psi\in \widetilde{\Phi}$ with $H(t)$ increasing, then the above inequalities hold with $\leq$ replacing by $\geq$.    

\noindent (ii) Assume that $\phi$ and $ \psi$ satisfy condition (d) $H(t)$ concave increasing with either $\phi, \psi\in \widetilde{\Phi}$ or $\phi, \psi\in \widetilde{\Psi}$. Then, for all $K\in \cS_0$,  \begin{eqnarray*}
\frac{\DOrliczA(K)}{n|K|}  \leq   H\bigg(\frac{\widetilde{\O}_{\psi}^{orlicz}(K)}{n|K|} \bigg),   \ \ \ \& \ \ \ \   \frac{\DOrliczG(K)}{n|K|} \leq   H\bigg(\frac{\widetilde{G}_{\psi}^{orlicz}(K)}{n|K|} \bigg). \end{eqnarray*} While if $\phi$ and $ \psi$ satisfy one of the following conditions:  (e) $H(t)$ convex decreasing with one in $\widetilde{\Phi}$ and another one in $\widetilde{\Psi}$; (f) $H(t)$  convex increasing with either $\phi, \psi\in \widetilde{\Phi}$ or $\phi, \psi\in \widetilde{\Psi}$,  then the above inequalities hold with $\leq$ replacing by $\geq$.  \et

 \noindent {\bf Proof.} We  only prove  the case $\DOrliczG(K)$ and omit the proof for $\DOrliczA(K)$.     \vskip 2mm \noindent (i).  For condition (a) $\phi\in \widetilde{\Phi}$ and $\psi\in \widetilde{\Psi}$ with $H(t)$ increasing and condition (b)  $\phi, \psi\in \widetilde{\Phi}$ with $H(t)$ decreasing: by Proposition \ref{bounded by volume}, one has, for $K\in \widetilde{\cK}$, 
\begin{eqnarray*} \frac{ \DOrliczG(K)}{ n|K|}  \leq \phi\left( {vrad(K^\circ)}\right) =H\left[\psi\left( {vrad(K^\circ)}\right)\right] \leq  H\bigg(\frac{ \widetilde{G}_{\psi}^{orlicz}(K)}{ n|K|} \bigg). \end{eqnarray*}    
 If functions $\phi\in \widetilde{\Psi}$ and $\psi\in \widetilde{\Phi}$ satisfy condition (c),  then  by Proposition \ref{bounded by volume},  
 the above inequalities hold with $\leq$ replacing by $\geq$.  
 
\vskip 2mm  \noindent (ii). 
For condition (d): the concavity of $H(t)$ with Jensen's inequality imply that, $\forall L\in \cS_0$,  \begin{eqnarray*}\frac{\dualorlicz(K, L)}{|K|} = \frac{1}{n|K|} \int _{S^{n-1}} H\bigg[\psi\!\left(\!\frac{\r_L(u)}{\r_K(u)}\!\right)\!\bigg] \r_K^n(u)\,d\s(u) \leq H\bigg(\frac{\widetilde{V}_{\psi}(K, L)}{|K|}\bigg).\end{eqnarray*} Let $H(t)$ be increasing and concave: by formula (\ref{Orlicz-geominimal-surface-area-convex}), one has, for $\phi,\psi\in \widetilde{\Phi}$ and for all $K\in \cS_0$,
\begin{eqnarray*}\frac{\DOrliczG(K)}{n|K|} &=&\inf_{L\in \widetilde{\cK}} \frac{n\dualorlicz(K, \vrad(L^\circ)L)}{n|K|} \\ &\leq&   H\bigg(\! \inf_{L\in \widetilde{\cK}} \frac{n\widetilde{V}_{\psi}(K, \vrad(L^\circ)L)}{n|K|}\!\bigg)=H\bigg(\frac{\widetilde{G}_{\psi}^{orlicz}(K)}{n|K|} \bigg). \end{eqnarray*} Replacing $\inf$ by $\sup$, one gets the analogous result for  $\phi,\psi\in \widetilde{\Psi}$,  due to formula  (\ref{Orlicz-geominimal-surface-area:concave}).

On the other hand, if $H(t)$ is convex, then  Jensen's inequality implies, \begin{eqnarray}\frac{\dualorlicz(K, \vrad(L^\circ)L)}{|K|}\geq  H\bigg(\frac{\widetilde{V}_{\psi}(K, \vrad(L^\circ)L)}{|K|}\bigg), \ \ \ \ \forall L\in \cS_0. \label{convex:decreasing---1} \end{eqnarray}  For $\phi \in \widetilde{\Psi}$ and $\psi\in \widetilde{\Phi}$ satisfy condition (e), i.e., $H(t)$ is convex and decreasing,   formulas (\ref{Orlicz-geominimal-surface-area-convex})-(\ref{Orlicz-geominimal-surface-area:concave}) imply that  $\forall K\in \cS_0$,  
 \begin{eqnarray*} \frac{\DOrliczG(K)}{n|K|}  & \geq&   \sup_{L\in \widetilde{\cK}}  H\bigg(\!\frac{n\widetilde{V}_{\psi}(K, \vrad(L^\circ)L)}{n|K|}\!\bigg) \\ &=& H
 \bigg(\!\inf_{L\in \widetilde{\cK}} \frac{n\widetilde{V}_{\psi}(K, \vrad(L^\circ)L)}{n|K|}\!\bigg) = H\bigg(\!\frac{\widetilde{G}_{\psi}^{orlicz}(K)}{n|K|}\! \bigg);\end{eqnarray*} 
By interchanging $\inf$ and $\sup$, one gets the analogous result for  $\phi \in \widetilde{\Phi}$ and $\psi\in  \widetilde{\Psi}$ with $H(t)$ being convex and decreasing.  

 For  $\phi, \psi\in \widetilde{\Phi}$ satisfying condition (f), i.e.,  $H(t)$ is convex increasing: by inequality (\ref{convex:decreasing---1}) and formula (\ref{Orlicz-geominimal-surface-area-convex}), one has, for all $K\in \cS_0$,  
 \begin{eqnarray*} \frac{\DOrliczG(K)}{n|K|}   &\geq&   \inf_{L\in \widetilde{\cK}}  H\bigg(\!\frac{n\widetilde{V}_{\psi}(K, \vrad(L^\circ)L)}{n|K|}\!\bigg)\\ &=&   H
 \bigg(\!\inf_{L\in \widetilde{\cK}} \frac{n\widetilde{V}_{\psi}(K, \vrad(L^\circ)L)}{n|K|}\!\bigg)=H\bigg(\!\frac{\widetilde{G}_{\psi}^{orlicz}(K)}{n|K|} \! \bigg).\end{eqnarray*}  
 Replacing $\inf$ by $\sup$, one gets the analogous result for $\phi, \psi\in \widetilde{\Psi}$ with $H(t)$ convex increasing. 
  
\bt\label{cyclic-pqr}
Let $q, r, s\neq n$ be such that either $s<r<0<q<n$,  or $0<s<r<q<n$, or $0<s<n<r<q$. Then,  for all $K\in \cS_0$, 
\begin{eqnarray*} 
\ \ \widetilde{G}^{orlicz}_r(K)  \leq  \big[\widetilde{G}^{orlicz}_q(K)\big] ^{\frac{r-s}{q-s}}
\big[\widetilde{G}^{orlicz}_s (K)\big]^{\frac{q-r}{q-s}}, \ \ \&\ \ \widetilde{\O}^{orlicz}_r(K)  \leq   \big[\widetilde{\O}^{orlicz}_q(K)\big] ^{\frac{r-s}{q-s}}
\big[\widetilde{\O}^{orlicz}_s (K)\big]^{\frac{q-r}{q-s}}.
\end{eqnarray*} 
\et

\vskip 2mm \noindent {\bf Proof.} We only prove the geominimal case and the affine case follows along the same line. Let $K\in  \cS_0$ and  $s<r<q$  (hence $0<\frac{q-r}{q-s}<1$).  H\"{o}lder's inequality (see \cite{HLP}) implies that    \begin{eqnarray}
n\widetilde{V}_r(K, Q)  &\leq & \big[n\widetilde{V}_s(K, Q)\big]^{\frac{q-r}{q-s}} \ \big[n\widetilde{V}_q(K, Q)\big]^{\frac{r-s}{q-s}},  \ \ \ \ \forall Q\in \cS_0. \label{mixed:p:surface:holder}
\end{eqnarray} 

\noindent Case (i). Let $s<r<0<q<n$.  Note that $t^q\in \widetilde{\Psi}$ as $0<q<n$. Then, for all $Q\in \widetilde{\cK}$, one has, 
  \begin{eqnarray*}
[\widetilde{G}^{orlicz}_q(K) ]^{\frac{r-s}{q-s}}\geq \big[ n \widetilde{V}_q(K, \vrad(Q^\circ)Q)\big]^{\frac{r-s}{q-s}}.\end{eqnarray*}  Note that $t^r, t^s\in \widetilde{\Phi}$ as $r, s<0$. Together with inequality (\ref{mixed:p:surface:holder}),  one has, 
\begin{eqnarray}
 \widetilde{G}^{orlicz}_r(K)&=&\inf_{Q\in \widetilde{\cK}} \big\{ n\widetilde{V}_r(K, \vrad(Q^\circ)Q) \big\} \nonumber\\ &\leq& [\widetilde{G}^{orlicz}_q(K) ] ^{\frac{r-s}{q-s}} \times \inf_{Q\in \widetilde{\cK}} \big\{ n\widetilde{V}_s(K, \vrad(Q^\circ)Q) \big\}^{\frac{q-r}{q-s}}\nonumber\\ &=& \big[\widetilde{G}^{orlicz}_q(K) \big] ^{\frac{r-s}{q-s}} \  \big[\widetilde{G}^{orlicz}_s(K) \big] ^{\frac{q-r}{q-s}}. \nonumber\end{eqnarray}  
The analogous result for the case  $0<s<n<r<q$ follows along the same line, if one notices that   $t^s\in \widetilde{\Psi}$  and  $t^q, t^r\in \widetilde{\Phi}$.

\vskip 2mm \noindent Case (ii). Let $0<s<r<q<n$, which clearly implies $t^q, t^r, t^s\in \widetilde{\Psi}$. Taking the supremum over $Q\in \widetilde{\cK}$ from both sides of inequality (\ref{mixed:p:surface:holder}), one gets,  for all $K\in \cS_0$, 
\begin{eqnarray*}
 \widetilde{G}^{orlicz}_r(K)   &\leq&   \sup_{Q\in \widetilde{\cK}} \big\{ n\widetilde{V}_q(K, \vrad(Q^\circ)Q) \big\}^{\frac{r-s}{q-s}} \sup_{Q\in \widetilde{\cK}} \big\{ n\widetilde{V}_s(K, \vrad(Q^\circ)Q) \big\}^{\frac{q-r}{q-s}} \\ &=&   \big[\widetilde{G}^{orlicz}_q(K) \big] ^{\frac{r-s}{q-s}} \  \big[\widetilde{G}^{orlicz}_s(K) \big] ^{\frac{q-r}{q-s}}.\end{eqnarray*}

  \section{Dual Orlicz mixed $L_{\phi}$ affine and geominimal surface areas } \label{section dual mixed}
  Various affine and geominimal surface areas for multiple convex bodies  have been studied extensively in, e.g., \cite{Lut1987, WernerYe2010, Y, Ye2014, DYzhuzhou2014}.  
In this section, the dual Orlicz mixed $L_{\phi}$ affine and geominimal surface areas for multiple star bodies are  briefly discussed.  Most of the proofs are either similar to those for single star body in Section \ref{subsection dual affine} or similar to those in \cite{Ye2014, DYzhuzhou2014}, and hence will be omitted.  
 
Let $\vec{\phi}=(\phi_1, \phi_2, \cdots, \phi_n)$ and $\vec{\phi}\in \widetilde{\Phi}^n$ (or $\vec{\phi}\in \widetilde{\Psi}^n$) means that each $\phi_i\in \widetilde{\Phi}$ (or $\phi_i\in \widetilde{\Psi}$). Similarly, $\cL=(L_1, \cdots, L_n)\in \cS_0^n$ means that each $L_i\in \cS_0$.  Define  $\widetilde{V}_{\vec{\phi}}(\bK, \cL)$ for $\bK, \cL\in \cS_0^n$ by \be\nonumber \widetilde{V}_{\vec{\phi}}(\bK, \cL)=\frac{1}{n}\int _{S^{n-1}}\prod_{i=1}^n \bigg[\phi_i\left(\frac{\r_{L_i}(u)}{\r_{K_i}(u)}\right)[\r_{K_i}(u)]^n \bigg]^{\frac{1}{n}}  \,d\s(u).\ee  When $\phi_i=\phi$, 
$K_i=K$ and $L_i=L$ for all $i=1, 2, \cdots, n$, one gets $ \widetilde{V}_{\vec{\phi}} (\bK;\cL)=\widetilde{V}_{\phi}(K, L). $  

We now propose definitions for the dual Orlicz mixed $L_{\phi}$ affine and geominimal surface areas. 
\bd\label{equivalent:mixed affine:surface:area-1} Let $K_1, \cdots,
K_n\in \cS_0$.   For $\vec{\phi}\in\widetilde{\Phi}^n$, define   $\DOrliczAmix(\bK)$ and $\DOrliczGmix(\bK)$ by
 \begin{eqnarray*} \DOrliczAmix(\bK)&=& \inf_{ \cL\in \widetilde{\cS}^n }
 \big\{n \widetilde{V}_{\vec{\phi}} (\bK;  \cL) \ \ with \ \  |L_1^\circ|=|L_2^\circ|=\cdots =|L_n^\circ|=\o_n\big\}, \\ \DOrliczGmix(\bK)&=& \inf_{ \cL\in \widetilde{\cK}^n }
 \big\{n \widetilde{V}_{\vec{\phi}}(\bK;  \cL) \ \ with \ \  |L_1^\circ|=|L_2^\circ|=\cdots =|L_n^\circ|=\o_n\big\}.  \end{eqnarray*} For   $\vec{\phi}\in \widetilde{\Psi}^n$,  $\DOrliczAmix(\bK)$ and   $\DOrliczGmix(\bK)$ are defined as above, but with $\inf$ replacing by $\sup$. 
\ed

   \noindent {\bf Remark.} As in \cite{DYzhuzhou2014}, one may be able to define several different dual Orlicz mixed $L_{\phi}$ affine and geominimal surface areas for $\bK$. In this paper, only the one defined by Definition  \ref{equivalent:mixed affine:surface:area-1} will be discussed and properties for others are very similar. Due to $\widetilde{\cK}^n\subset \widetilde{\cS}^n$, for $\bK\in \cS_0^n$, one has, $\DOrliczAmix(\bK)\leq \DOrliczGmix(\bK)$ for $\vec{\phi}\in \widetilde{\Phi}^n$  and $\DOrliczAmix(\bK)\geq \DOrliczGmix(\bK)$ for  $\vec{\phi}\in  \widetilde{\Psi}^n$. Moreover, the dual Orlicz mixed $L_{\phi}$ affine and geominimal surface areas are affine invariant: for $\bK\in \cS_0^n$ and for $\vec{\phi}\in  \widetilde{\Phi}^n\cup  \widetilde{\Psi}^n$,   $$ \DOrliczAmix(T\bK)=\DOrliczAmix(\bK); \ \ \ \DOrliczGmix(T\bK)=\DOrliczGmix(\bK), \ \  \forall T\in SL(n),$$  where  $T \bK=(T K_1, \cdots, T K_n)$ for $T\in SL(n)$.  For $\bK\in \cS_0^n$ and $\vec{\phi}\in  \widetilde{\Phi}^n$, one has \begin{eqnarray*}  
\big[\DOrliczAmix(\bK)\big]^{n} \leq  \big[\DOrliczGmix(\bK)\big]^{n} \leq  \widetilde{S}_{\phi_1}(K_1)\cdots \widetilde{S}_{\phi_n}(K_n). \end{eqnarray*}

\bt \label{Alexander-Fenchel:both} Let $\bK\in \cS_0^n $. For $\vec{\phi}\in  \widetilde{\Phi}^n\cup  \widetilde{\Psi}^n$, one has \begin{eqnarray*}  
\big[\DOrliczAmix(\bK)\big]^{n} \leq \prod_{i=1}^n \widetilde{\O}_{\phi_i}^{orlicz}(K_i) \ \ and \ \ \  \big[\DOrliczGmix(\bK)\big]^{n}  \leq \prod_{i=1}^n  \widetilde{G}_{\phi_i}^{orlicz}(K_i). \end{eqnarray*} Moreover, if $\vec{\phi}\in  \widetilde{\Psi}^n$, the following Alexander-Fenchel type inequalities hold:  Let $m$ be an integer such that  $1 \leq m \leq n$,  then
\begin{eqnarray*}
\big[\DOrliczAmix(\bK)\big]^{m}&\leq& \prod_{i=0}^{m-1}\widetilde{\O}_{(\phi_{1},
\cdots, \phi_{n-m}, \phi_{n-i},\cdots, \phi_{n-i})}^{orlicz}(K_{1},
\cdots, K_{n-m}, \underbrace{K_{n-i},\cdots, K_{n-i}}_{m}),  \\  \big[\DOrliczGmix(\bK)\big]^{m}&\leq& \prod_{i=0}^{m-1}\widetilde{G}_{(\phi_{1},
\cdots, \phi_{n-m}, \phi_{n-i},\cdots, \phi_{n-i})}^{orlicz}(K_{1},
\cdots, K_{n-m}, \underbrace{K_{n-i},\cdots, K_{n-i}}_{m}).
\end{eqnarray*}  \et

 \noindent \textbf{Proof.} We only prove the geominimal case and omit the proof for the affine case.  In fact,  H\"{o}lder's inequality (see \cite{HLP}) implies 
\begin{eqnarray}
 \big[\widetilde{V}_{\vec{\phi}}(\bK;  \cL)\big]^m\!\!\!\! &\leq& \!\!\!\! \frac{1}{n} \prod_{i=0}^{m-1} \int_{S^{n-1}}  \bigg[\phi_{n-i}\left(\frac{\r_{L_{n-i}}(u)}{\r_{K_{n-i}}(u)}\right) [\r_{K_{n-i}}(u)]^n  \bigg]^{\frac{m}{n}}   \prod_{j=1}^{n-m} \bigg[\phi_j\left(\frac{\r_{L_j}(u)}{\r_{K_j}(u)}\right) [\r_{K_j}(u)]^n \bigg]^{\frac{1}{n}} \,d\s(u) \nonumber \\ \!\!\!\!&=&\!\!\!\! \prod_{i=0}^{m-1}\widetilde{V}_{(\phi_{1}, \cdots , \phi_{n-m}, \phi_{n-i},\cdots, \phi_{n-i}) } (K_{1},\!\cdots\!,\! K_{n-m}, \underbrace{K_{n-i},\!\cdots\!, K_{n-i}}_{m}; L_{1},\!\cdots\!, L_{n-m}, \underbrace{L_{n-i},\!\cdots\!, L_{n-i}}_{m}). \nonumber \end{eqnarray}  Let $\vec{\phi}\in \widetilde{\Psi}^n$. 
 Taking the supremum over $\cL\in \widetilde{\cK}^n$ with $|L_1^\circ|=\cdots =|L_n^{\circ}|=\o_n$, one gets the desired Alexander-Fenchel type inequality if one notices that for all $\cL\in \widetilde{\cK}^n$ and all $i=0, \cdots, m-1$,  \begin{eqnarray*}  n\widetilde{V}_{(\phi_{1}, \cdots , \phi_{n-m}, \phi_{n-i},\cdots, \phi_{n-i}) } (K_{1},\!\cdots\!,\! K_{n-m}, \underbrace{K_{n-i},\!\cdots\!, K_{n-i}}_{m}; L_{1},\!\cdots\!, L_{n-m}, \underbrace{L_{n-i},\!\cdots\!, L_{n-i}}_{m})\\ \leq \widetilde{G}_{(\phi_{1},
\cdots, \phi_{n-m}, \phi_{n-i},\cdots, \phi_{n-i})}^{orlicz}(K_{1},\!
\cdots\!, K_{n-m}, \underbrace{K_{n-i},\!\cdots\!, K_{n-i}}_{m}).  \end{eqnarray*}  
   Note that if $m=n$, then  $\big[\widetilde{V}_{\vec{\phi}}(\bK;  \cL)\big]^n \leq  \prod_{i=1}^{n} \widetilde{V}_{\phi_i} (K_{i},  L_{i}).$ 
 Definitions  \ref{Orlicz affine surface} and \ref{equivalent:mixed affine:surface:area-1} imply that for $\vec{\phi}\in  \widetilde{\Phi}^n$\begin{eqnarray*}
\big[\DOrliczGmix(\bK)\big]^n &=& \bigg[\inf_{ \cL \in \widetilde{\cK}^n}
 \big\{n \widetilde{V}_{\vec{\phi}}(\bK;  \cL) \ \ with \ \  |L_1^\circ|=|L_2^\circ|=\cdots =|L_n^\circ|=\o_n\big\}\bigg]^n \\  &\leq &  \prod_{i=1}^{n} \inf_{L_i \in \widetilde{\cK}} \big\{ n \widetilde{V}_{\phi_i} (K_{i},  L_{i}) \ \ with \ \  |L_i^\circ| =\o_n\big\} =\prod_{i=1}^n \widetilde{G}_{\phi_i}^{orlicz}(K_i).\end{eqnarray*}  
 Replacing $\inf$ by $\sup$, one gets the desired result for  $\phi\in  \widetilde{\Psi}^n$.  
 
 The following Orlicz affine isoperimetric type inequalities follows from  Theorems \ref{isoperimetric:geominimal} and \ref{Alexander-Fenchel:both}.  

\bt \label{isoperimetric:geominimal:mixed case} 
 Let $\bK\in \cS_0^n$. 
 
 \vskip 2mm \noindent (i) 
 For $\vec{\phi}\in \widetilde{\Phi}_1^n$, one has 
\begin{eqnarray*}
\big[\DOrliczAmix(\bK)\big]^{n} &\leq&    \prod_{i=1}^n \widetilde{\O}^{orlicz}_{\phi_i}\big([B_{(K_i)^\circ}]^\circ\big), \ \ \ \ \bK\in \widetilde{\cS} ^n;  \\ \big[\DOrliczAmix(\bK)\big]^{n} &\leq&  \big[\DOrliczGmix(\bK)\big]^{n} \leq \prod_{i=1}^n \widetilde{G}^{orlicz}_{\phi_i}([B_{(K_i)^\circ}]^\circ), \ \ \ \ \bK\in \widetilde{\cK}^n.  
\end{eqnarray*}  

\noindent (ii)  For $\vec{\phi}\in \widetilde{\Psi}^n$ and $\bK\in \cS_0^n$, one has\begin{eqnarray*}
\big[\DOrliczGmix(\bK)\big]^{n} \leq  \big[\DOrliczAmix(\bK)\big]^{n} \leq   \prod_{i=1}^n \widetilde{\O}^{orlicz}_{\phi_i}(B_{K_i})= \prod_{i=1}^n \widetilde{G}^{orlicz}_{\phi_i}(B_{K_i}).
\end{eqnarray*}   
 \et 
 
For $\bK\in \cS_0^n$, write $\widetilde{G}^{orlicz}_p(\bK)$ for $\DOrliczGmix(\bK)$ and $\widetilde{\O}^{orlicz}_p(\bK)$ for $\DOrliczAmix(\bK)$ if $\vec{\phi}=(t^p, \cdots, t^p)$. Similar to the proof of  Theorem \ref{cyclic-pqr}, one has the following theorem. 
 \bt
Let $q, r, s\neq n$ be such that either $s<r<0<q<n$,  or $0<s<r<q<n$, or $0<s<n<r<q$. Then,  for $\bK\in \cS_0^n$, 
\begin{eqnarray*} 
  \ \ \widetilde{G}^{orlicz}_r(\bK)  \leq  \big[\widetilde{G}^{orlicz}_q(\bK)\big] ^{\frac{r-s}{q-s}}
\big[\widetilde{G}^{orlicz}_s (\bK)\big]^{\frac{q-r}{q-s}}, \ \  \&\ \  \widetilde{\O}^{orlicz}_r(\bK)  \leq   \big[\widetilde{\O}^{orlicz}_q(\bK)\big] ^{\frac{r-s}{q-s}}
\big[\widetilde{\O}^{orlicz}_s (\bK)\big]^{\frac{q-r}{q-s}}.
\end{eqnarray*} 
\et

 For $i\in \bbR$, define $\widetilde{V}_{\phi_1, \phi_2, i}(K, L; Q_1, Q_2)$ with $K, L, Q_1, Q_2\in \cS_0$  by  
   \be \nonumber n \widetilde{V}_{\phi_1, \phi_2, i}(K, L; Q_1, Q_2) \!=\!\! \int _{S^{n-1}}\!\! \bigg[\phi_1\!\!\left(\!\frac{\r_{Q_1}(u)}{\r_{K}(u)}\!\right)\! [\r_{K}(u)]^n\bigg]^{\frac{n-i}{n}}\!\bigg[\phi_2\!\!\left(\!\frac{\r_{Q_2}(u)}{\r_{L}(u)}\!\right)\! [\r_{L}(u)]^n\bigg]^{\frac{i}{n}} \!\! \,d\s(u).\ee 
 The dual Orlicz $i$-th mixed $L_{\phi}$ affine and geominimal surface areas for $K, L\in \cS_0$, denoted by $\widetilde{\O}_{\phi_1, \phi_2,i}^{orlicz} (K, L) $ and $\widetilde{G}_{\phi_1, \phi_2,i}^{orlicz} (K, L)$ respectively,  are defined as follows.  \bd\label{mixed i th affine surface}  Let $K, L\in \cS_0$ and $i\in \bbR$.      For $\phi_1, \phi_2\in \widetilde{\Phi}$, we define  \begin{eqnarray*} \widetilde{\O}_{\phi_1, \phi_2,i}^{orlicz} (K, L) 
 &=&\inf_{ \{Q_1, Q_2\in \widetilde{\cS} \} }
 \big\{n\widetilde{V}_{\phi_1, \phi_2,i}\big(K, L; Q_1, Q_2\big): \ \ |Q_1^\circ|=|Q_2^\circ|=\o_n \big\}, \\ \widetilde{G}_{\phi_1, \phi_2,i}^{orlicz} (K, L) 
&=&\inf_{ \{Q_1, Q_2\in \widetilde{\cK}\} }
 \big\{n\widetilde{V}_{\phi_1, \phi_2,i}\big(K, L; Q_1, Q_2 \big): \ \ |Q_1^\circ|=|Q_2^\circ|=\o_n \big\}.
  \end{eqnarray*} While if $\phi_1, \phi_2\in \widetilde{\Psi}$,   $\widetilde{\O}_{\phi_1, \phi_2,i}^{orlicz} (K, L) $ and respectively $\widetilde{G}_{\phi_1, \phi_2,i}^{orlicz} (K, L)$  are defined as above, but with $\inf$ replacing by $\sup$.   \ed    
 The dual Orlicz $i$-th mixed $L_{\phi}$ affine and  geominimal surface areas  are all affine invariant. Moreover, for $K, L\in\cS_0$  and $i\in \bbR$, one has, due to $\widetilde{\cK}\subset \widetilde{\cS}$, \begin{eqnarray} \widetilde{\O}_{\phi_1, \phi_2,i}^{orlicz} (K, L)&\leq& \widetilde{G}_{\phi_1, \phi_2, i}^{orlicz} (K, L), \ \ \phi_1, \phi_2\in \widetilde{\Phi};\label{compare:ith:----1}\\ \widetilde{\O}_{\phi_1, \phi_2,i}^{orlicz} (K, L)&\geq& \widetilde{G}_{\phi_1, \phi_2, i}^{orlicz} (K, L), \ \ \phi_1, \phi_2\in \widetilde{\Psi}.\label{compare:ith:----2}  \end{eqnarray}
 
\bt \label{cyclic-theorem} Let $K, L\in\cS_0$ and $i<j<k$. For $\phi_1, \phi_2\in\widetilde{\Psi}$, one has 
 \begin{eqnarray*} \big[\widetilde{\O}_{\phi_1, \phi_2,j}^{orlicz} (K, L)\big] ^{k-i}&\leq& \big[\widetilde{\O}_{\phi_1, \phi_2,i}^{orlicz} (K, L)\big]^{k-j}\big[\widetilde{\O}_{\phi_1, \phi_2,k}^{orlicz} (K, L)\big]^{j-i}; \\
  \big[\widetilde{G}_{\phi_1, \phi_2,j}^{orlicz} (K, L)\big] ^{k-i}&\leq& \big[\widetilde{G}_{\phi_1, \phi_2,i}^{orlicz} (K, L)\big]^{k-j}\big[\widetilde{G}_{\phi_1, \phi_2,k}^{orlicz} (K, L)\big]^{j-i}. \end{eqnarray*}    \et

 \noindent \textbf{Proof.}   Let $i< j<k$ which implies
 $0<\frac{k-j}{k-i}<1$. H\"{o}lder's inequality  implies that,
 \begin{eqnarray}\nonumber
\widetilde{V}_{\phi_1, \phi_2, j}(K, L; Q_1, Q_2)\leq [\widetilde{V}_{\phi_1, \phi_2, i}(K, L; Q_1, Q_2)]^{\frac{k-j}{k-i}}[\widetilde{V}_{\phi_1, \phi_2, k}(K, L; Q_1, Q_2)]^{\frac{j-i}{k-i}}.   
\end{eqnarray}
The desired result follows by taking the supremum over $Q_1, Q_2\in \widetilde{\cS}$ and $Q_1, Q_2\in \widetilde{\cK}$ respectively with $|Q_1^\circ|=|Q_2^\circ|=\o_n$. 
   
 \bt Let $K, L\in \cS_0$.  \vskip 2mm \noindent (i)  Let $0\leq i \leq n$ and $\phi_1, \phi_2\in \widetilde{\Phi}_1$. One has  \begin{eqnarray*}\big[\widetilde{\O} ^{orlicz}_{\phi_1, \phi_2, i}(K, L)\big]^{n}&\leq& \big[\widetilde{G} ^{orlicz}_{\phi_1, \phi_2, i}(K, L)\big]^{n} \leq  \big[\widetilde{G}_{\phi_1}^{orlicz}([B_{K^\circ}]^\circ)\big]^{n-i}\big[\widetilde{G}_{\phi_2}^{orlicz}([B_{L^\circ}]^\circ)\big]^{i}, \ \ \ K\in \widetilde{\cK}; \\ \big[\widetilde{\O} ^{orlicz}_{\phi_1, \phi_2, i}(K, L)\big]^{n}& \leq&  \big[\widetilde{\O}_{\phi_1}^{orlicz}([B_{K^\circ}]^\circ)\big]^{n-i}\big[\widetilde{\O}_{\phi_2}^{orlicz}([B_{L^\circ}]^\circ)\big]^{i}, \ \ \ K\in \widetilde{\cS}.
   \end{eqnarray*}  
 
\noindent (ii)  Let $0\leq i \leq n$ and $\phi_1, \phi_2\in \widetilde{\Psi}$. One has, for $K, L\in \cS_0$,  \begin{eqnarray*}\big[\widetilde{G} ^{orlicz}_{\phi_1, \phi_2, i}(K, L)\big]^{n}\leq \big[\widetilde{\O} ^{orlicz}_{\phi_1, \phi_2, i}(K, L)\big]^{n} \leq  \big[\widetilde{\O}_{\phi_1}^{orlicz}(B_K)\big]^{n-i}\big[\widetilde{\O} _{\phi_2}^{orlicz}(B_L)\big]^{i}.
 \end{eqnarray*} (iii) Let $\E$ be an origin-symmetric ellipsoid and $\phi_1, \phi_2\in \widetilde{\Psi}$. For $i>n$ and $K\in \cS_0$, one has,    
\begin{eqnarray*} \big[\widetilde{\O}_{\phi_1, \phi_2,i}^{orlicz} (K, \E)\big]^n \geq \big[\widetilde{G}_{\phi_1, \phi_2, i}^{orlicz} (K, \E) \big]^n \geq   \big[\widetilde{G}_{\phi_1}^{orlicz}(B_K) \big]^{n-i}\big[\widetilde{G} _{\phi_2}^{orlicz}(\E)\big]^{i}. 
 \end{eqnarray*}  \et

 \noindent  {\bf Proof.}  By H\"{o}lder's inequality (see
\cite{HLP})  and Definitions \ref{Orlicz affine surface} and \ref{mixed i th affine surface}, one has, for $0\leq i\leq n$, $K, L\in \cS_0$  and $\phi_1, \phi_2\in \widetilde{\Phi}$ or $\phi_1, \phi_2\in\widetilde{\Psi}$,  \begin{eqnarray}
 \big[\widetilde{\O} ^{orlicz}_{\phi_1, \phi_2, i}(K, L)\big]^{n} &\leq&  \big[\widetilde{\O} _{\phi_1}^{orlicz}(K)\big]^{n-i}\big[\widetilde{\O} _{\phi_2}^{orlicz}(L)\big]^{i},  \label{dual mixed volume-100} \\
\big[\widetilde{G} ^{orlicz}_{\phi_1, \phi_2, i}(K, L)\big]^{n}&\leq& \big[\widetilde{G} _{\phi_1}^{orlicz}(K)\big]^{n-i}\big[\widetilde{G} _{\phi_2}^{orlicz}(L)\big]^{i}. \label{dual mixed volume-1}
 \end{eqnarray}
Similarly, for  $\phi_1, \phi_2\in\widetilde{\Psi}$,   $K, L\in \cS_0$ and  $i< 0$ or $i>n$,  one has, \begin{eqnarray}
\big[\widetilde{\O} ^{orlicz}_{\phi_1, \phi_2, i}(K, L)\big]^{n}&\geq& \big[\widetilde{\O} _{\phi_1}^{orlicz}(K)\big]^{n-i}\big[\widetilde{\O} _{\phi_2}^{orlicz}(L)\big]^{i}, \nonumber \\ \big[\widetilde{G} ^{orlicz}_{\phi_1, \phi_2, i}(K, L)\big]^{n}&\geq& \big[\widetilde{G} _{\phi_1}^{orlicz}(K)\big]^{n-i}\big[\widetilde{G} _{\phi_2}^{orlicz}(L)\big]^{i}.\label{Geominimal:ball:ith---1}
\end{eqnarray}

 \vskip 2mm \noindent (i). Let $\phi_1, \phi_2\in \widetilde{\Phi}_1$  and $0\leq i\leq n$. Combining  inequality (\ref{dual mixed volume-1})  with Theorem \ref{isoperimetric:geominimal} and inequality (\ref{compare:ith:----1}),  one gets, for $K, L\in \widetilde{\cK}$,   \begin{eqnarray*}\big[\widetilde{\O} ^{orlicz}_{\phi_1, \phi_2, i}(K, L)\big]^{n} & \leq & \big[\widetilde{G} ^{orlicz}_{\phi_1, \phi_2, i}(K, L)\big]^{n}\leq  \big[\widetilde{G} _{\phi_1}^{orlicz}(K)\big]^{n-i}\big[\widetilde{G}  _{\phi_2}^{orlicz}(L)\big]^{i}  \\ &\leq& \big[\widetilde{G} _{\phi_1}^{orlicz}([B_{K^\circ}]^\circ)\big]^{n-i}\big[\widetilde{G} _{\phi_2}^{orlicz}([B_{L^\circ}]^\circ)\big]^{i}.
  \end{eqnarray*}   Similarly, combining  inequality (\ref{dual mixed volume-100})  with Theorem \ref{isoperimetric:geominimal},  one gets, for $K, L\in \widetilde{\cS}$, \begin{eqnarray*} \big[\widetilde{\O} ^{orlicz}_{\phi_1, \phi_2, i}(K, L)\big]^{n}& \leq&  \big[\widetilde{\O}_{\phi_1}^{orlicz}([B_{K^\circ}]^\circ)\big]^{n-i}\big[\widetilde{\O}_{\phi_2}^{orlicz}([B_{L^\circ}]^\circ)\big]^{i}, \ \ \ K\in \widetilde{\cS}.
   \end{eqnarray*} 

\noindent (ii).  Let $\phi_1, \phi_2\in\widetilde{\Psi}$  and $0\leq i\leq n$. Combining  inequality (\ref{dual mixed volume-100})  with Theorem \ref{isoperimetric:geominimal} and inequality (\ref{compare:ith:----2}),  one gets, for $K, L\in \widetilde{\cS}$,   \begin{eqnarray*}\big[\widetilde{G} ^{orlicz}_{\phi_1, \phi_2, i}(K, L)\big]^{n} & \leq & \big[\widetilde{\O} ^{orlicz}_{\phi_1, \phi_2, i}(K, L)\big]^{n}\leq  \big[\widetilde{\O} _{\phi_1}^{orlicz}(K)\big]^{n-i}\big[\widetilde{\O}  _{\phi_2}^{orlicz}(L)\big]^{i}  \\ &\leq& \big[\widetilde{\O} _{\phi_1}^{orlicz}(B_K)\big]^{n-i}\big[\widetilde{\O} _{\phi_2}^{orlicz}(B_L)\big]^{i}.
  \end{eqnarray*}      
  \noindent (iii). Let $i>n$ and $\phi_1, \phi_2\in \widetilde{\Psi}$.  
 Inequalities (\ref{compare:ith:----2}) and (\ref{Geominimal:ball:ith---1}) together with Theorem \ref{isoperimetric:geominimal} imply  \begin{eqnarray*} \big[\widetilde{\O}_{\phi_1, \phi_2,i}^{orlicz} (K, \E)\big]^n &\geq& \big[\widetilde{G}_{\phi_1, \phi_2, i}^{orlicz} (K, \E) \big]^n   \geq  \big[\widetilde{G} _{\phi_1}^{orlicz}(K)\big]^{n-i}\big[\widetilde{G} _{\phi_2}^{orlicz}(\E)\big]^{i} \\ & \geq&  \big[\widetilde{G}_{\phi_1}^{orlicz}(B_K) \big]^{n-i}\big[\widetilde{G} _{\phi_2}^{orlicz}(\E)\big]^{i}. 
    \end{eqnarray*}

\vskip 2mm \noindent {\bf Acknowledgments.} The research of DY is supported
 by a NSERC grant.

\vskip 2mm \noindent Deping Ye, \ \ \ {\small \tt deping.ye@mun.ca}\\
{\small \em Department of Mathematics and Statistics\\
   Memorial University of Newfoundland\\
   St. John's, Newfoundland, Canada A1C 5S7 }

\end{document}